\documentclass[10pt, a4paper, twoside]{article}
\usepackage[portuges]{babel}
\usepackage[latin1]{inputenc}

\usepackage{amsmath,amsfonts}
\usepackage{amsthm}

\newcommand{\newsec}[1]{\section{#1}\setcounter{equation}{0}
                                    \setcounter{thmTEMA}{0}
                                    \setcounter{teoTEMA}{0}
                                    \setcounter{lemmaTEMA}{0}
                                    \setcounter{lemaTEMA}{0}
                                    \setcounter{defTEMAi}{0}
                                    \setcounter{defTEMAp}{0}
                                    \setcounter{propTEMAi}{0}
                                    \setcounter{propTEMAp}{0} }

\newtheorem{teoTEMA}{Teorema}

\newtheorem{defTEMAp}{Defini\c{c}\~{a}o}

\newtheorem{coroTEMAp}{Corol\'ario}

\renewenvironment{abstract}
{\begin{list}{}{%
  \setlength{\rightmargin}{0.5cm}
  \setlength{\leftmargin}{0.5cm}
  \small}
  \item[] }
{\end{list}}

\newcommand{\runningheads}{\markboth}
\normalsize

\newtheorem{obsTEMAp}{Observação}

\theoremstyle{definition}
\newtheorem{explTEMAp}{Exemplo}

\begin{document}


\title{Computação Algébrica no Cálculo das Variações:\\
    Determinação de Simetrias e Leis de Conservação\footnote{Partially
    presented at XXVII CNMAC (Congresso Nacional de Matemática Aplicada
    e Computacional), FAMAT/PUCRS, Porto Alegre, RS, Brasil, 13-16 September
    2004. Supported by the program PRODEP III/5.3/2003 and the R\&D unit CEOC.}}

\author{Paulo~D.~F.~Gouveia\\
     \texttt{pgouveia@ipb.pt}
     \and
     Delfim~F.~M.~Torres\\
     \texttt{delfim@mat.ua.pt}}

\date{Control Theory Group (\textsf{cotg})\\
     Centre for Research in Optimization and Control (\textsf{CEOC})\\
     University of Aveiro, Department of Mathematics\\
     3810-193 Aveiro, Portugal}

\maketitle

\runningheads
{Gouveia e Torres} {Computação Algébrica no Cálculo das Variações}

\setlength\arraycolsep{1pt}


\begin{abstract}
{\bf Resumo}. Os problemas de optimização dinâmica (em espaços de funções)
tratados pelo cálculo das variações, são normalmente resolvidos
por recurso às condições necessárias de Euler-Lagrange,
que são equações diferenciais de segunda ordem
(ou de ordem superior, quando os problemas variacionais envolvem
derivadas de ordem superior a um). Estas equações são, em geral,
não lineares e de difícil resolução. Uma forma de as simplificar
consiste em obter leis de conservação:
primeiros integrais das equações diferenciais de Euler-Lagrange.
Os primeiros integrais permitem baixar a ordem das equações
e, em casos extremos, com um número suficientemente grande
de primeiros integrais independentes, resolver o problema por
completo. Se em áreas como a Física e a Economia
a questão da existência de leis de conservação é
resolvida de forma bastante natural, a própria aplicação
sugerindo as leis de conservação (\textrm{e.g.} conservação
de energia, conservação da quantidade de movimento, conservação
do rendimento, etc.), de um ponto de vista estritamente matemático,
dado um problema do cálculo de variações, o processo
de obtenção das leis de conservação ou, até mesmo,
a demonstração de que elas existem (ou não), deixa de ser uma questão óbvia.
Neste trabalho mostramos como um sistema de computação algébrica
como o Maple pode ser muito útil na abordagem a estas questões.
Apresentamos um conjunto de facilidades computacionais
simbólicas que permitem, de uma forma sistemática e automática,
identificar as leis de conservação de uma dada funcional integral
do cálculo das variações. O algoritmo usado tem por base o célebre
teorema de Emmy Noether, que associa a existência de leis de conservação
às propriedades de invariância do problema
(à existência de simetrias variacionais).
Vários exemplos ilustrativos são apresentados,
mostrando a utilidade das ferramentas desenvolvidas.
\end{abstract}


\begin{abstract}
{\bf Abstract}. The dynamic optimization problems treated
by the calculus of variations are usually solved with
the help of the 2nd order Euler-Lagrange differential equations.
These equations are, generally speaking, nonlinear, and very hard to solve.
One way to address the problem is to obtain conservation laws of lower order
than those of the corresponding Euler-Lagrange equations.
While in Physics and Economics the question of existence of conservation laws
is treated in a rather natural way, because the application itself suggest
the conservation laws (e.g., conservation of energy, income/health law),
from a strictly mathematical point of view, given a problem
of the calculus of variations, it is not obvious how one might
derive a conservation law or, for that matter, if it even has
a conservation law. The question we address is thus to develop
computational facilities, based on a systematic method,
which permits to identify functionals that have conservation laws.
The central result we use is the celebrated Noether's theorem.
This theorem links conservation laws with the invariance properties of the problem
(with symmetries), and provides an algorithm for finding conservation laws.
Thus the problem is reduced to the one of finding the variational symmetries.
We show how a Computer Algebra System can help to find the symmetries
and the conservation laws in the calculus of variations.
Several illustrative examples are given.
\end{abstract}


\smallskip



\newsec{Introdução}
\setcounter{obsTEMAp}{0}

O cálculo das variações é uma área clássica da
matemática, com mais de três séculos de idade,
extremamente activa no século XXI, e com inúmeras
aplicações práticas na mecânica, economia, ciências
dos materiais, ciências do espaço e engenharia.
O cálculo das variações está na origem de muitas áreas
mais recentes, como sejam a análise funcional
e o controlo óptimo \cite{EDP,CVFAwOC,Leitao}.

A computação algébrica, também chamada de
computação científica ou simbóli\-ca,
permite trabalhar com expressões matemáticas de
maneira simbólica, não numé\-rica, e é uma
área de investigação moderna, que surgiu na segunda
metade do século XX. Se é verdade que os actuais sistemas
de computação algébrica são extremamente poderosos, também não é
menos verdade que muito existe por fazer no sentido
de se colocarem os computadores a realizar tarefas
matemáticas verdadeiramente interessantes.

O uso da teoria do cálculo das variações,
na resolução de problemas concretos,
requer cálculos não numéricos:
no cálculo das variações a presença de fenómenos como
o de Lavrentiev \cite{fenomLav}, faz com que
as soluções algébricas exactas sejam muito
mais convenientes do que as soluções numéricas.
Tendo em conta que a resolução dos problemas do cálculo
das variações passa quase obrigatoriamente pela resolução
das equações diferenciais de Euler-Lagrange,
equações estas, em geral, de difícil resolução,
neste trabalho propomos um conjunto de procedimentos
computacionais algébricos que permitem automatizar
o processo de obtenção de leis de conservação
(primeiros integrais das equações de Euler-Lagrange).
Como é bem conhecido da teoria das equações diferenciais,
estas leis de conservação são de extrema utilidade,
permitindo baixar a ordem das equações \cite{OLV}.

O presente trabalho encontra-se organizado em sete secções.
Em \S\ref{sec:CV} fazemos uma breve apresentação
do cálculo das variações, da condição necessária de Euler-Lagrange
(definida em Maple na \S\ref{sec:ProcMaple})
e definimos \emph{lei de conservação}.
Em \S\ref{sec:SV} definimos \emph{simetria}
(\emph{invariância} de um problema do cálculo das variações)
e apresentamos um método construtivo para as obter
(implementado em Maple na \S\ref{sec:ProcMaple},
por intermédio do procedimento \emph{Simetria}).
Em \S\ref{sec:LC} formulamos o célebre princípio de Noether,
que nos dá uma fórmula explícita para as leis de conservação
do cálculo das variações, em função
dos geradores que definem a simetria. O teorema de Noether
é definido em Maple na \S\ref{sec:ProcMaple} por intermédio
do procedimento \emph{Noether}.
Vários exemplos, mostrando o modo de uso e a utilidade das funções Maple
definidas em \S\ref{sec:ProcMaple}, são apresentados
em \S\ref{sec:EI}. Por fim, apresentamos em \S\ref{sec:conclusao}
algumas conclusões e perspectivas de trabalho futuro.


\newsec{Cálculo das Variações}
\label{sec:CV}
\setcounter{obsTEMAp}{0}

Neste trabalho consideramos problemas do cálculo das variações de
ordem superior: minimizar uma funcional integral
\begin{equation}
\label{eq:funcional}
J[\mathbf{x}(\cdot)] = \int_{a}^{b}
L(t,\mathbf{x}(t),\dot{\mathbf{x}}(t),\ldots,\mathbf{x}^{(m)}(t))
\,\textrm{d}t \, ,
\end{equation}
onde o Lagrangeano $L$ é uma função real que assumimos ser continuamente
diferenciável em $[a, b] \times\mathbb{R}^{n \times (m+1)}$; $t \in
\mathbb{R}$, é a variável independente; $\mathbf{x}(t)=[ x_1(t)\ x_2(t)$
$\cdots x_n(t) ]^\textrm{T} \in \mathbb{R}^n$, as variáveis dependentes;
$\mathbf{x}^{(i)}(t)=[ \frac{\mathrm{d} ^i x_1(t)}{\mathrm{d} t^i} \
\frac{\mathrm{d}^ix_2(t)}{\mathrm{d}t^i} \cdots
\frac{\mathrm{d}^ix_n(t)}{\mathrm{d}t^i}  ]^\textrm{T} \in \mathbb{R}^n$, com
$i=1,\ldots, m$, as derivadas de ordem $i$ das variáveis dependentes em ordem
a $t$; e $\dot{\mathbf{x}}(t)\equiv \mathbf{x}^{(1)}(t)$.
Assim, o Lagrangeano considerado, depende explicitamente de uma variável
independente, de $n$ variáveis dependentes e das suas $m$ primeiras derivadas.
Para $m = 1$ obtemos o problema fundamental do cálculo das variações.

Na minimização da funcional (\ref{eq:funcional}) é usual recorrer-se
ao sistema de equações de Euler-Lagrange
\begin{equation}
\label{eq:EL}
\frac{\partial{L}}{\partial{\mathbf{x}}}
+ \sum_{i=1}^{m}(-1)^i\frac{\textrm{d}^i}
{\textrm{d}t^i} \left(\frac{\partial{L}}{\partial{\mathbf{x}^{(i)}}} \right)
= \mathbf{0},
\end{equation}
onde $\frac{\partial L}{\partial \mathbf{x}^{(i)}}
=[ \frac{\partial L}{\partial x_1^{(i)}}
\ \frac{\partial L}{\partial x_2^{(i)}}
\ \cdots \ \frac{\partial L}{\partial x_n^{(i)}}]$
e $L$, e suas derivadas,
são avaliadas ao longo de $(t,\mathbf{x}(t),
\dot{\mathbf{x}}(t),\ldots,\mathbf{x}^{(m)}(t))$.

\begin{defTEMAp}
Às soluções da equação de Euler-Lagrange (\ref{eq:EL})
chamamos \emph{extremais}.
\end{defTEMAp}

\begin{obsTEMAp}
Para o problema fundamental do cálculo das variações,
\textrm{i.e.}, quando a funcional não depende de derivadas
de $\mathbf{x}(t)$ de maior ordem do que a primeira ($m=1$),
o sistema de equações de Euler-Lagrange (\ref{eq:EL}) reduz-se a
\begin{equation*}
\frac{\partial{L}}{\partial{\mathbf{x}}} -\frac{\emph{d}}
{\emph{d}t} \left(\frac{\partial{L}}{\partial{\dot{\mathbf{x}}}} \right)
= \mathbf{0} \, .
\end{equation*}
\end{obsTEMAp}

\begin{obsTEMAp}
Se a funcional não envolver mais do que uma variável dependente
($n=1$), todas as grandezas presentes em (\ref{eq:funcional}) e
(\ref{eq:EL}) são escalares.
\end{obsTEMAp}

Em \S\ref{sec:ProcMaple} definimos o procedimento
\emph{EulerLagrange} que tem por entrada o Lagrangeano
e como saída o sistema de equações de Euler-Lagrange correspondente,
que resulta da aplicação da equação (\ref{eq:EL}).

\begin{defTEMAp}
Uma função $t \rightarrow
\phi(t,\mathbf{x}(t),\dot{\mathbf{x}}(t),\ldots,\mathbf{x}^{(k)}(t))$,
$k < 2m$, que se mantenha constante ao longo de todas as extremais
do problema (\ref{eq:funcional}), é chamada de \emph{primeiro integral}
(de ordem $k$) da equação de Euler-Lagrange (\ref{eq:EL}).
À equação
\begin{equation*}
\phi(t,\mathbf{x}(t),\dot{\mathbf{x}}(t),\ldots,\mathbf{x}^{(k)}(t))
= \text{const}
\end{equation*}
chamamos \emph{lei de conservação} (de ordem $k$).
\end{defTEMAp}

As leis de conservação são muito úteis, pois permitem
reduzir a ordem das equações diferenciais
de Euler-Lagrange. Com um número suficientemente elevado
de primeiros integrais independentes,
é mesmo possível determinar explicitamente as extremais.
Consideremos, a título de exemplo,
o seguinte problema ($n = m = 1$):
\begin{equation*}
J[x(\cdot)] = \int_{a}^{b} \left(\dot{x}^2(t) - x^2(t)\right)
\,\textrm{d}t \longrightarrow \min \, .
\end{equation*}
Neste caso o Lagrangeano é dado por $L(x,\dot{x}) = -x^2 + \dot{x}^2$
e a equação de Euler-Lagrange \eqref{eq:EL}
reduz-se a $\ddot{x}(t) + x(t) = 0$.
Se multiplicarmos esta equação diferencial por $\cos(t)$
obtemos $\frac{d}{dt} \left(\dot{x} \cos(t) + x \sin(t)\right) = 0$;
enquanto que se a multiplicarmos por $-\sin(t)$ obtemos
$\frac{d}{dt} \left(-\dot{x} \sin(t) + x \cos(t)\right) = 0$.
Temos então duas leis de conservação:
\begin{equation}
\label{eq:extrivialLC}
\begin{cases}
\dot{x} \cos(t) + x \sin(t) &= c_1 \, ,\\
-\dot{x} \sin(t) + x \cos(t) &= c_2 \, .
\end{cases}
\end{equation}
Resulta de imediato das duas leis de conservação \eqref{eq:extrivialLC}
que as extremais têm a forma
$x(t) = c_1 \sin(t) + c_2 \cos(t)$. Claro que este exemplo é trivial:
a equação diferencial de Euler-Lagrange pode facilmente ser resolvida
sem recurso a leis de conservação, uma vez que ela é linear.
Usando o Maple e o nosso procedimento  \emph{EulerLagrange} faríamos:
\small
\begin{verbatim}
> L := v^2 - x^2:
> dsolve(EulerLagrange(L,t,x,v));
\end{verbatim}
\begin{equation*}
 \left\{ x \left( t \right) ={\it \_C1}\,\sin \left( t \right)
 +{\it \_C2}\,\cos \left( t \right)  \right\}
\end{equation*}
\normalsize
No caso não linear o comando Maple \emph{dsolve} nem sempre é capaz
de obter soluções explícitas (\textrm{e.g.} para
a equação não linear de Emden-Fowler considerada no Exemplo~\ref{ex:Emden-Fowler};
ou para o problema de Kepler tratado no Exemplo~\ref{ex:Kepler}).
Nessas situações as leis de conservação são muito úteis.
Coloca-se então a seguinte questão: \emph{Dada uma funcional integral
do tipo (\ref{eq:funcional}), como obter leis de conservação?}
Esta questão foi resolvida por Emmy Noether em 1918 \cite{Noether}:
se a funcional for invariante sob determinado tipo de transformações
(transformações de invariância ou simetrias), então existem fórmulas explícitas
para as leis de conservação. O Teorema de Noether encontra muitas
aplicações em campos concretos da Engenharia, como a Sismologia
e a Metalurgia \cite{alik}. A dificuldade na sua aplicação reside
na obtenção das simetrias (na obtenção das transformações de invariância).
Neste trabalho automatizamos, por recurso ao sistema de computação algébrica
Maple \cite{Lenimar}, a determinação de simetrias e a correspondente aplicação
do Teorema de Noether.


\newsec{Simetrias Variacionais}
\setcounter{obsTEMAp}{0}
\label{sec:SV}

Para o estudo das propriedades de invariância das
funcionais do cálculo das variações considera-se
uma família uni-paramétrica de transformações
$\mathbf{h}^s(t,\mathbf{x})$ que forma um grupo local de Lie
\cite{LOG,BRU}. A família uni-paramétrica $\mathbf{h}^s(t,\mathbf{x})$
representa um conjunto de $n+1$ transformações
de $[a, b] \times \mathbb{R}^n$ em  $\mathbb{R} \times \mathbb{R}^n$
\begin{equation}
\label{eq:transf}
t^s = h_t^s(t,\mathbf{x})\textrm{,}\quad
x_i^s = h^s_{x_i}(t,\mathbf{x})\textrm{, com}\;\, i=1,\ldots,n,
\end{equation}
a que correspondem $n+1$ geradores infinitesimais representados por
\begin{equation}
T(t,\mathbf{x})  =  \left.
\frac{\partial}{\partial{s}} h^s_t(t,\mathbf{x})\right|_{s=0}\textrm{,}\quad
X_i(t,\mathbf{x})  =  \left.
\frac{\partial}{\partial{s}} h^s_{x_i}(t,\mathbf{x})\right|_{s=0}\textrm{, com}
\, i=1,\ldots,n. \label{eq:geradores}
\end{equation}

\begin{defTEMAp}
\label{defin:invar}
A funcional (\ref{eq:funcional}) diz-se invariante
no intervalo $[a,b]$ sob as transformações uni-paramétricas
(\ref{eq:transf}) se, para todo o $s$ suficientemente pequeno,
\[
\int_{\alpha}^{\beta}{L(t,\mathbf{x}(t),\dot{\mathbf{x}}(t),
\ldots,\mathbf{x}^{(m)}(t))\,\emph{d}t}
= \int_{\alpha^s}^{\beta^s}{L(t^s,\mathbf{x}^s(t^s),\dot{\mathbf{x}}^s(t^s),
\ldots,{\mathbf{x}^s}^{(m)}(t^s))\,\emph{d}t^s},
\]
em qualquer subintervalo $[\alpha,\beta] \subseteq [a,b]$;
com $\alpha^s=h_t^s(\alpha,\mathbf{x}(\alpha))$
e $\beta^s=h_t^s(\beta,\mathbf{x}(\beta))$.
\end{defTEMAp}

\begin{obsTEMAp}
Nas condições da definição \ref{defin:invar}
as trans\-for\-ma\-ções uni-pa\-ra\-mé\-tri\-cas
(\ref{eq:transf}) constituem uma \emph{simetria variacional}
da funcional (\ref{eq:funcional}).
\end{obsTEMAp}

O teorema que se segue estabelece uma condição necessária
e suficiente de invariância, de extrema importância
para os objectivos a que nos propomos.

\begin{teoTEMA}[\cite{TOR}]
\label{teor:gerad}
A funcional (\ref{eq:funcional}) é invariante
sob as transformações uni-paramétricas (\ref{eq:transf}),
com geradores infinitesimais $T$ e
$\mathbf{X}$ (\ref{eq:geradores}), se, e apenas se,
\begin{equation}
\label{eq:gerad}
\frac{\partial L}{\partial t}T+ \sum_{i=0}^{m}
\frac{\partial L}{\partial \mathbf{x}^{(i)}} \cdot
\mathbf{p}^i+ L\frac{\emph{d}T}{\emph{d}t}=0,
\end{equation}
onde
\begin{equation}
\label{eq:gerad1}
\mathbf{p}^0=\mathbf{X}\textrm{,} \quad
\mathbf{p}^{i+1}=\frac{\emph{d}\mathbf{p}^{i}}{\emph{d}t}
-\mathbf{x}^{(i+1)}\frac{\emph{d}T}{\emph{d}t}
\textrm{,}\;\, i=0,\ldots,m-1 \, .
\end{equation}
Em (\ref{eq:gerad}) e (\ref{eq:gerad1}) assumimos
que $T$ e $\mathbf{X}=\left[ X_1\ X_2
\cdots X_n \right]^T$ são avaliadas em função de $(t,\mathbf{x})$
e $\mathbf{p}^i$ em função de $(t,\mathbf{x},
\dot{\mathbf{x}},\ldots,\mathbf{x}^{(i)})$, com $i=0,1,\ldots,m$.
\end{teoTEMA}

\begin{coroTEMAp}
Quando o Lagrangeano $L$ não depende de derivadas de $\mathbf{x}(t)$
de ordem superior à primeira ($m=1$), a equação (\ref{eq:gerad})
toma a forma
\[
\frac{\partial L}{\partial t}T+\frac{\partial L}{\partial \mathbf{x}}
\cdot \mathbf{X}+\frac{\partial L}{\partial \dot{\mathbf{x}}}
\cdot \left( \frac{\emph{d}\mathbf{X}}{\emph{d}t}
-\dot{\mathbf{x}}\frac{\emph{d}T}{\emph{d}t} \right) +
L\frac{\emph{d}T}{\emph{d}t}=0,
\]
com
\[
\frac{\textrm{d}T}{\textrm{d}t} = \frac{\partial T}{\partial t}
+ \frac{\partial T}{\partial
\mathbf{x}}\cdot\dot{\mathbf{x}}\textrm{,}\quad
\frac{\emph{d}\mathbf{X}}{\emph{d}t}
= \frac{\partial \mathbf{X}}{\partial t} + \frac{\partial \mathbf{X}}
{\partial \mathbf{x}}\cdot\dot{\mathbf{x}},
\]
onde
\[
\frac{\partial \mathbf{X}}{\partial \mathbf{x}}  =
\left[\frac{\partial \mathbf{X}}{\partial x_1}\
\frac{\partial \mathbf{X}}{\partial x_2}\
\cdots\  \frac{\partial \mathbf{X}}{\partial x_n}\right]
= \left[
\begin{array}{cccc}
\frac{\partial X_1}{\partial x_1} & \frac{\partial X_1}{\partial x_2}
& \cdots & \frac{\partial X_1}{\partial x_n}\\
\frac{\partial X_2}{\partial x_1} & \frac{\partial X_2}{\partial x_2}
& \cdots & \frac{\partial X_2}{\partial x_n}\\
\vdots&\vdots&\ddots&\vdots\\
\frac{\partial X_n}{\partial x_1} & \frac{\partial X_n}{\partial x_2}
& \cdots & \frac{\partial X_n}{\partial x_n}
\end{array}
\right].
\]
\end{coroTEMAp}

\begin{obsTEMAp}
Todas as derivadas totais presentes em (\ref{eq:gerad}) e (\ref{eq:gerad1})
podem ser expressas por derivadas parciais,
usando as igualdades $\frac{\textrm{d}T}{\textrm{d}t}
= \frac{\partial T}{\partial t} + \frac{\partial T}{\partial \mathbf{x}}\cdot\dot{\mathbf{x}}$ e
\begin{eqnarray*}
\frac{\textrm{d}\mathbf{p}^{i}}{\textrm{d}t} = \frac{\partial
\mathbf{p}^{i}}{\partial t} + \sum_{k=0}^{i}\frac{\partial
\mathbf{p}^{i}}{\partial \mathbf{x}^{(k)}}\cdot\mathbf{x}^{(k+1)} \, , \quad
i = 0,\ldots,m-1 \, ,
\end{eqnarray*}
onde
\[
\frac{\partial \mathbf{p}^{i}}{\partial \mathbf{x}^{(k)}}  =
\left[\frac{\partial \mathbf{p}^{i}}{\partial x_1^{(k)}}
\ \frac{\partial \mathbf{p}^{i}}{\partial x_2^{(k)}}\ \cdots\
\frac{\partial \mathbf{p}^{i}}{\partial x_n^{(k)}}\right]
= \left[
\begin{array}{cccc}
\frac{\partial p_1^{i}}{\partial x_1^{(k)}}
& \frac{\partial p_1^{i}}{\partial x_2^{(k)}}
& \cdots & \frac{\partial p_1^{i}}{\partial x_n^{(k)}}\\
\frac{\partial p_2^{i}}{\partial x_1^{(k)}}
& \frac{\partial p_2^{i}}{\partial x_2^{(k)}}
& \cdots & \frac{\partial p_2^{i}}{\partial x_n^{(k)}}\\
\vdots&\vdots&\ddots&\vdots\\
\frac{\partial p_n^{i}}{\partial x_1^{(k)}}
& \frac{\partial p_n^{i}}{\partial x_2^{(k)}}
& \cdots & \frac{\partial p_n^{i}}{\partial x_n^{(k)}}
\end{array}
\right].
\]
\end{obsTEMAp}

O teorema \ref{teor:gerad}, para além de servir de teste
à existência de simetrias, estabelece um algoritmo
para a determinação dos correspondentes
geradores infinitesimais. Como veremos, este facto é crucial:
o teorema de Noether (Teorema~\ref{teor:lc}) afirma que
as leis de conservação associadas a uma dada simetria variacional
apenas dependem dos geradores infinitesimais.

Dado um Lagrangeano $L$, determinamos os geradores infinitesimais
$T$ e $\mathbf{X}$ de uma família uni-paramétrica de
transformações simétricas pelo seguinte método.
A equação (\ref{eq:gerad}) é uma equação diferencial
nas $n+1$ funções incógnitas $T$, $X_1$, $X_2$, \dots, e $X_n$,
que pretendemos determinar. Porém, a equação
tem de permanecer válida para todos os $x_i$, $i=1,\ldots,n$.
Como as funções $T$, $X_1$, $X_2$, \dots, e $X_n$
dependem de $t$ e $x_i$, $i=1,\ldots,n$,
ao substituirmos, na equação (\ref{eq:gerad}),
$L$ e todas as suas derivadas parciais pelos seus valores,
obtemos um polinómio nas $n\times m$ variáveis
$\dot{x}_1,\ldots, \dot{x}_n$, $x_1^{(2)},\ldots,
x_n^{(2)},\ldots$,$x_1^{(m)},\ldots, x_n^{(m)}$
e suas potências. Para que a equação seja válida para
todos os valores das variáveis do polinómio,
todos os seus coeficientes devem ser nulos.
Notamos que os termos do polinómio poderão ser em maior
número que as incógnitas do problema ($n+1$), pelo que
a condição necessária e suficiente (\ref{eq:gerad})
pode conduzir a um sistema de equações sem solução.
Tal facto significa apenas que nem todas as funcionais
integrais do cálculo das variações admitem simetrias variacionais
(ver Exemplo~\ref{ex:Thomas-Fermi}). O sistema de equações a resolver,
para a obtenção dos geradores, é um sistema de equações
diferenciais às derivadas parciais. No entanto,
ao contrário das equações diferenciais ordinárias
de Euler-Lagrange, em geral não lineares e de difícil resolução,
este sistema é linear em $\frac{\partial T}{\partial t}$,
$\frac{\partial T}{\partial \mathbf{x}}$, $\frac{\partial \mathbf{X}}{\partial t}$ e
$\frac{\partial \mathbf{X}}{\partial \mathbf{x}}$.

A resolução do sistema de equações diferencias parciais que deriva
da expressão (\ref{eq:gerad}), nomeadamente quando lidamos
com valores de $n$ e $m$ superiores à unidade,
envolve um número muito elevado de cálculos, o que torna premente
a necessidade de nos munirmos de ferramentas computacionais
que automatizem o trabalho. Com esse fim,
desenvolvemos um procedimento em \texttt{Maple},
designado \emph{Simetria}. O respectivo código, assim como
o dos restantes procedimentos por nós desenvolvidos,
é apresentado na secção~\ref{sec:ProcMaple}.
O procedimento \emph{Simetria}
tem por entrada a expressão que caracteriza o Lagrangeano
e como saída os respectivos geradores infinitesimais.
No caso do Lagrangeano não admitir
qualquer simetria, obtemos geradores nulos
(\textrm{cf.} Exemplo~\ref{ex:Thomas-Fermi}).

Na secção que se segue mostramos como os geradores,
obtidos por intermédio do nosso procedimento
\emph{Simetria}, podem ser usados na obtenção
explícita de leis de conservação.


\newsec{Leis de Conservação}
\label{sec:LC}
\setcounter{obsTEMAp}{0}

Emmy Noether foi a primeira a estabelecer a relação
entre a existência de simetrias e a
existência de leis de conservação. Esta ligação
constitui um princípio universal, passível
de ser formulado na forma de teorema
nos mais diversos contextos e sob as mais variadas
hipóteses \cite{alik,LOG77,OLV,delfimEJC,delfimCAO,TOR}.

\begin{teoTEMA}[Teorema de Noether \cite{TOR}]
\label{teor:lc}
Se a funcional (\ref{eq:funcional}) é invariante
sob as transformações uni-paramétricas (\ref{eq:transf}),
com geradores infinitesimais $T$ e $\mathbf{X}$, então
\begin{equation}
\label{eq:lc}
\sum_{i=1}^m\Psi^i\cdot \mathbf{p}^{i-1} + \left(L-\sum_{i=1}^m\Psi^i
\cdot \mathbf{x}^{(i)} \right)T
= \,\textrm{const,} \quad t \in [a,b],
\end{equation}
com
\begin{eqnarray*}
\Psi^m&=&\frac{\partial L}{\partial \mathbf{x}^{(m)}},\nonumber\\
\Psi^{i-1}&=&\frac{\partial L}{\partial \mathbf{x}^{(i-1)}}
-\frac{\emph{d}\Psi^{i}}{\emph{d}t}\textrm{,}\quad i=m,m-1,\ldots,2\quad\quad\\
\frac{\emph{d}\Psi^{i}}{\emph{d}t}&=&\frac{\partial \Psi^{i}}{\partial t}
+ \sum_{k=0}^{2m-i}\left(\mathbf{x}^{(k+1)}\right)^\emph{T}\cdot\frac{\partial
\Psi^{i}}{\partial \mathbf{x}^{(k)}},
\end{eqnarray*}
onde
\begin{eqnarray*}
\frac{\partial \Psi^i}{\partial \mathbf{x}^{(k)}}
=\left[
\begin{array}{c}
\frac{\partial \Psi^i}{\partial x_1^{(k)}}\\
\frac{\partial \Psi^i}{\partial x_2^{(k)}}\\
\vdots\\ \frac{\partial \Psi^i}{\partial x_n^{(k)}}
\end{array}
\right]
=
\left[
\begin{array}{cccc}
\frac{\partial \psi_1^i}{\partial x_1^{(k)}}
& \frac{\partial \psi_2^i}{\partial x_1^{(k)}}
& \cdots & \frac{\partial \psi_n^i}{\partial x_1^{(k)}}\\
\frac{\partial \psi_1^i}{\partial x_2^{(k)}}
& \frac{\partial \psi_2^i}{\partial x_2^{(k)}} & \cdots
& \frac{\partial \psi_n^i}{\partial x_2^{(k)}}\\
\vdots&\vdots&\ddots&\vdots\\
\frac{\partial \psi_1^i}{\partial x_n^{(k)}}
& \frac{\partial \psi_2^i}{\partial x_n^{(k)}} & \cdots
& \frac{\partial \psi_n^i}{\partial x_n^{(k)}}
\end{array}
\right]
\end{eqnarray*}
e onde se considera a grandeza $\Psi^i$ avaliada em
$(t,\mathbf{x}(t),\dot{\mathbf{x}}(t),\ldots, \mathbf{x}^{(2m-i)}(t))$,
$i=1,\ldots,m$.
\end{teoTEMA}

\begin{coroTEMAp}
\label{cor:TeoNoether}
Quando o Lagrangeano $L$ não depende das derivadas de $\mathbf{x}(t)$
de maior ordem do que a primeira ($m=1$),
a equação (\ref{eq:lc}) reduz-se a
\begin{equation*}
\frac{\partial L}{\partial \dot{\mathbf{x}}} \cdot \mathbf{X}
+ \left(L - \frac{\partial L}{\partial \dot{\mathbf{x}}} \cdot
\dot{\mathbf{x}}\right)T=\,\textrm{const.}
\end{equation*}
\end{coroTEMAp}

As leis de conservação que procuramos são obtidas
substituindo nas equações (\ref{eq:gerad1}) e (\ref{eq:lc}) os
geradores infinitesimais $T$ e $\mathbf{X}$
encontrados pelo método descrito na secção anterior.
Na secção~\ref{sec:ProcMaple} definimos
o procedimento \emph{Noether}. Este procedimento
tem por entradas o Lagrangeano e os geradores infinitesimais,
que são obtidos por intermédio do nosso procedimento
\emph{Simetria}, e como saída a correspondente
lei de conservação \eqref{eq:lc}.
Resumindo: dado um problema do cálculo das variações \eqref{eq:funcional},
obtemos as leis de conservação, de uma forma automática,
através de um processo de duas etapas:
com o nosso procedimento \emph{Simetria} obtemos
todas as possíveis simetrias do problema;
recorrendo depois ao nosso procedimento \emph{Noether},
baseado no teorema~\ref{teor:lc}, obtemos
as correspondentes leis de conservação.
Na secção seguinte apresentamos alguns exemplos
que ilustram todo o processo.


\newsec{Exemplos Ilustrativos}
\label{sec:EI}
\setcounter{obsTEMAp}{0}

Consideramos agora várias situações concretas,
mostrando a funcionalidade e a utilidade
das ferramentas desenvolvidas.

\begin{explTEMAp}
\label{ex:Ex1}
Começamos com um exemplo muito simples em que o Lagrangeano
depende apenas de uma variável dependente ($n = 1$) e não existem
derivadas de ordem superior à primeira ($m = 1$): $L(t,x,\dot{x})=t\dot{x}^2$.

\noindent Com a definição \texttt{Maple}
\small
\begin{verbatim}
> L:= t*v^2;
\end{verbatim}
\[
L := tv^2\nonumber
\]
\normalsize
o nosso procedimento \emph{Simetria}
determina os geradores infinitesimais
das simetrias do problema
do cálculo das variações em consideração:
\small
\begin{verbatim}
> Simetria(L,t,x,v);
\end{verbatim}
\[
\{ T \left( t,x \right) = \left( 2\,{\it \_C1}\,\ln  \left( t \right)
+{\it \_C3} \right) t, \;
X \left( t,x \right) ={\it \_C1}\,x+{\it \_C2}\}
\]
\normalsize
A família de geradores depende de três parâmetros que advêm
das constantes de integração. Com as substituições\footnote{O sinal
de percentagem (\%) é um operador usado em \texttt{Maple}
para referenciar o resultado do comando anterior.
Para uma introdução ao \texttt{Maple}, veja-se \cite{Lenimar}.}

\small
\begin{verbatim}
> subs(_C1=1,_C2=0,_C3=0,%);
\end{verbatim}
\[
\left\{ T \left( t,x \right) =2\,
\ln\left( t \right) t,X \left( t,x \right) =x \right\}
\]
\normalsize
obtemos os geradores descritos em \cite[pp.~210 e 214]{BRU}.
A lei de conservação correspondente a estes geradores
é facilmente obtida por intermédio do nosso procedimento \emph{Noether}:
\small
\begin{verbatim}
> LC := Noether(L,t,x,v,%);
\end{verbatim}
\[
LC :=
x \left( t \right) t{\frac {d}{dt}}x \left( t \right)
-{t}^{2} \left( {\frac {d}{dt}}x \left( t \right)\right)^{2}
\ln\left( t \right) ={\it const}
\]
\normalsize
É, neste caso, muito fácil verificar a validade
da lei de conservação obtida. Por definição,
basta mostrar que a igualdade é verificada ao longo
das extremais. A equação de Euler-Lagrange
é a equação diferencial de $2^a$ ordem
\small
\begin{verbatim}
> EulerLagrange(L,t,x,v);
\end{verbatim}
\[
\left\{ -2\,{\frac {d}{dt}}x \left( t \right)
-2\,t{\frac {d^2}{dt^2}  }x \left( t \right) =0 \right\}
\]
\normalsize
e as extremais são as suas soluções:
\small
\begin{verbatim}
> dsolve(%);
\end{verbatim}
\begin{equation*}
\{x(t) = {\it \_C1}+{\it \_C2} \ln(t)\}
\end{equation*}
\normalsize
Substituindo as extremais na lei de conservação, obtemos,
como esperado, uma proposição verdadeira:
\small
\begin{verbatim}
> expand(subs(%,LC));
\end{verbatim}
\begin{equation*}
{\it \_C2} {\it \_C1} = const
\end{equation*}
\normalsize
\end{explTEMAp}

\begin{explTEMAp}[Problema de Kepler]
\label{ex:Kepler}
Analisamos agora as simetrias e leis
de conservação do problema de Kepler \cite[p.~217]{BRU}.
Neste caso o Lagrangeano tem
duas variáveis dependentes ($n = 2$) e
não envolve derivadas de ordem superior ($m = 1$):
\begin{equation*}
L(t,\mathbf{q},\dot{\mathbf{q}})=\frac{m}{2}\left(\dot{q}_1^2
+\dot{q}_2^2\right)+\frac{K}{\sqrt{q_1^2+q_2^2}}.
\end{equation*}
Vamos determinar a fórmula geral das leis de conservação.
Neste caso não é possível validar a lei de conservação
por aplicação directa da definição, como fizemos para
o exemplo anterior, pois o Maple não é capaz
de resolver o respectivo sistema de equações
de Euler-Lagrange
\small
\begin{verbatim}
> L:=m/2*(v[1]^2+v[2]^2)+K/sqrt(q[1]^2+q[2]^2);
\end{verbatim}
\[
L\, := \,1/2\,m \left( {v_{{1}}}^{2}+{v_{{2}}}^{2} \right)
+{\frac {K}{\sqrt {{q_{{1}}}^{2}+{q_{{2}}}^{2}}}}
\]
\begin{verbatim}
> EulerLagrange(L,t,[q[1],q[2]],[v[1],v[2]]);
\end{verbatim}
\[
\left\{ -m{\frac {d^2}{d t^2 }}q_{{1}} \left( t \right)
-{\frac {Kq_{{1}} \left( t \right) }{ \left(  q_{{1}} \left( t \right)^{2}
+  q_{{2}} \left( t \right)  ^{2} \right) ^{3/2}}}=0,
\right. \;
\left. -m{\frac {d^2}{d t^2 }}q_{{2}} \left( t \right)
-{\frac {Kq_{{2}} \left( t \right) }{ \left(  q_{{1}} \left( t \right)^{2}
+  q_{{2}} \left( t \right) ^{2} \right) ^{3/2}}}=0 \right\}
\]
\begin{verbatim}
> Simetria(L, t, [q[1],q[2]], [v[1],v[2]]);
\end{verbatim}
\[
 \{ X_{{1}} \left( t,q_{{1}},q_{{2}} \right)
 ={\it \_C2}\,q_{{2}}, \; T \left( t,q_{{1}},q_{{2}} \right)
 ={\it \_C1},\;
X_{{2}} \left( t,q_{{1}},q_{{2}} \right)
=-{\it \_C2}\,q_{{1}} \}
\]
\begin{verbatim}
> Noether(L, t, [q[1],q[2]], [v[1],v[2]], %):
> expand(%);
\end{verbatim}
\begin{eqnarray}
{\it \_C2}\,q_{{2}} \left( t \right) m{\frac {d}{dt}}q_{{1}}
\left( t \right) -{\it \_C2}\,q_{{1}} \left( t \right) m{\frac {d}{dt}}q_{{2}}
\left( t \right)
-1/2\,{\it \_C1}\, \left( {\frac {d}{dt}}q_{{1}}
\left( t \right)\right)^{2}m\nonumber \\
-1/2\,{\it \_C1}\,
\left( {\frac {d}{dt}}q_{{2}} \left( t \right)\right)^{2}m
+{\frac {{\it \_C1}\,K}{\sqrt { q_{{1}} \left( t \right)^{2}
+ q_{{2}} \left( t \right) ^{2}}}}={\it const}
\nonumber
\end{eqnarray}
\normalsize
\end{explTEMAp}

\begin{explTEMAp}
Vejamos o caso de um Lagrangeano
com duas variáveis dependentes ($n = 2$)
e com derivadas de ordem superior ($m = 2$):
\begin{equation*}
L(t,\mathbf{x},\dot{\mathbf{x}},\ddot{\mathbf{x}})=\dot{x}_1^2+\ddot{x}_2^2
\end{equation*}
\small
\begin{verbatim}
> L:=v[1]^2+a[2]^2;
\end{verbatim}
\[
L\, := \,{v_{{1}}}^{2}+{a_{{2}}}^{2}
\]
\begin{verbatim}
> Simetria(L, t, [x[1],x[2]], [v[1],v[2]], [a[1],a[2]]);
\end{verbatim}
\begin{eqnarray}
&&\left\{\right.T \left( t,x_{{1}},x_{{2}} \right)
={\it \_C1}\,t+{\it \_C2}, \nonumber \\
&&X_{{1}} \left( t,x_{{1}},x_{{2}} \right)
=1/2\,{\it \_C1}\,x_{{1}}+{\it \_C5}, \;
X_{{2}} \left( t,x_{{1}},x_{{2}} \right)
=3/2\,{\it \_C1}\,x_{{2}}+{\it \_C3}\,t+{\it \_C4}\left.\right\}
\nonumber
\end{eqnarray}
\begin{verbatim}
> LC := Noether(L, t, [x[1],x[2]], [v[1],v[2]], [a[1],a[2]], %);
\end{verbatim}
\begin{eqnarray}
LC :=
&&2 \left( 1/2\,{\it \_C1}\,x_{{1}} \left( t \right)
+{\it \_C5} \right) {\frac {d}{dt}}x_{{1}} \left( t \right)
-2\, \left( 3/2\,{\it \_C1}\,x_{{2}} \left( t \right)
+{\it \_C3}\,t+{\it \_C4} \right) {\frac
{d^3}{d t^3 }}x_{{2}} \left( t \right)\nonumber\\
&+&2 \left( {\it \_C3}
+1/2\,{\it \_C1}\,{\frac{d}{dt}}x_{{2}}\left(t\right)\right)
{\frac {d^2}{d t^2  }}x_{{2}} \left( t \right)\nonumber\\
&+&
\bigg(- \left( {\frac {d}{dt}}x_{{1}} \left( t \right)\right)^{2}
- \left( {\frac {d^2}{d t^2}}x_{{2}} \left( t \right)  \right) ^{2}
+2 {\frac {d}{dt}}x_{{2}} \left( t \right)
{\frac {d^3}{d t^3}}x_{{2}} \left( t \right)  \bigg)
\left( {\it \_C1}\,t+{\it \_C2} \right)
={\it const}
\nonumber
\end{eqnarray}
\normalsize
Tal como para o Exemplo~\ref{ex:Ex1}, também aqui
é fácil verificar, por aplicação directa da definição,
a validade da lei de conservação obtida:
\small
\begin{verbatim}
> EulerLagrange(L,t, [x[1],x[2]], [v[1],v[2]], [a[1],a[2]]);
\end{verbatim}
\begin{displaymath}
\left\{ -2\,{\frac {d^2}{d t^2}}x_{{1}} \left( t \right)
=0,2\,{\frac {d^4}{d t^4}}x_{{2}} \left( t \right) =0 \right\}
\end{displaymath}
\begin{verbatim}
> dsolve(%);
\end{verbatim}
\[
 \left\{\right. x_{{1}} \left( t \right)
 ={\it \_C1}\,t+{\it \_C2},\;
 x_{{2}} \left( t \right) =1/6\,{\it \_C3}\,{t}^{3}
 +1/2\,{\it \_C4}\,{t}^{2}+{\it \_C5}\,t+{\it \_C6} \left.\right\}
\]
\begin{verbatim}
> expand(subs(%,LC));
\end{verbatim}
\[
2\,{\it \_C1}\,{\it \_C5}-3\,{\it \_C3}\,{\it \_C1}\,{\it \_C6}
+{\it \_C1}\,{\it \_C5}\,{\it \_C4}
-{{\it \_C4}}^{2}{\it \_C2}+2\,{\it \_C3}
\,{\it \_C5}\,{\it \_C2}={\it const}
\]
\normalsize
\end{explTEMAp}

\begin{explTEMAp}[Emden-Fowler]
\label{ex:Emden-Fowler}
Consideremos o problema variacional definido pelo Lagrangeano
\small
\begin{verbatim}
> L:= t^2/2*(v^2-(1/3)*x^6);
\end{verbatim}
\begin{equation*}
L:=\frac{{t}^{2} \left( {v}^{2}-\frac{{x}^{6}}{3} \right)}{2}
\end{equation*}
\normalsize
A respectiva equação diferencial de Euler-Lagrange é conhecida
na astrofísica como a equação de Emden-Fowler \cite[p. 220]{BRU}:
\small
\begin{verbatim}
> EL := EulerLagrange(L,t,x,v);
\end{verbatim}
\begin{equation*}
EL := \left\{ -2\,t{\frac {d}{dt}}x \left( t \right) -{t}^{2}{\frac {d^{2}}
{d{t}^{2}}}x \left( t \right) -{t}^{2} \left( x \left( t \right)
 \right) ^{5}=0 \right\}
\end{equation*}
\normalsize
Encontramos os geradores infinitesimais,
que conduzem a uma simetria variacional
para a funcional de Emden-Fowler,
por intermédio da nossa função \emph{Simetria}:
\small
\begin{verbatim}
> S := Simetria(L,t,x,v);
\end{verbatim}
\begin{equation*}
S := \left\{ X \left( t,x \right) =-\frac{x{\it \_C1}}{2},T \left( t,x \right)
={\it \_C1}\,t \right\}
\end{equation*}
\normalsize
Por exemplo,
\small
\begin{verbatim}
> S2 := subs(_C1=-6,S);
\end{verbatim}
\begin{equation*}
                S2 := \{T(t, x) = -6 t, X(t, x) = 3 x\}
\end{equation*}
\normalsize
Aplicando o Teorema de Noether (Teorema~\ref{teor:lc}),
estabelecemos a seguinte lei de conservação:
\small
\begin{verbatim}
> simplify(Noether(L,t,x,v,S2));
\end{verbatim}
\begin{equation*}
{t}^{2} \left( 3\,x \left( t \right) {\frac {d}{dt}}x \left( t
 \right) +3\, \left( {\frac {d}{dt}}x \left( t \right)  \right) ^{2}t+
t \left( x \left( t \right)  \right) ^{6} \right) ={\it const}
\end{equation*}
\normalsize
\end{explTEMAp}

\begin{explTEMAp}[Thomas-Fermi]
\label{ex:Thomas-Fermi}
Mostramos agora um exemplo de um problema do cálculo das variações
que não possui nenhuma simetria variacional. Seja
\small
\begin{verbatim}
> L:=1/2 * v^2 + 2/5 * (x^(5/2))/(sqrt(t));
\end{verbatim}
\begin{equation*}
L:=\frac{v^{2}}{2}+\frac{2\,x^{\frac{5}{2}}}{5\sqrt {t}}
\end{equation*}
\normalsize
A equação de Euler-Lagrange associada a este Lagrangeano corresponde
à equação diferencial de Thomas-Fermi \cite[p. 220]{BRU}:
\small
\begin{verbatim}
> EL := EulerLagrange(L,t,x,v);
\end{verbatim}
\begin{equation*}
EL :=  \left\{ -{\frac {d^{2}}{d{t}^{2}}}x \left( t \right) +{\frac {
 \left( x \left( t \right)  \right) ^{\frac{3}{2}}}{\sqrt {t}}}=0 \right\}
\end{equation*}
\normalsize
A nossa função \emph{Simetria} devolve, neste caso, geradores nulos.
Como explicado na \S\ref{sec:SV}, isto significa que este problema
do cálculo das variações não admite simetrias.
\small
\begin{verbatim}
> Simetria(L, t, x, v);
\end{verbatim}
\begin{equation*}
\{X(t,x) = 0, T(t,x) = 0\}
\end{equation*}
\normalsize
A função \emph{Noether} resulta num truísmo:
\small
\begin{verbatim}
> Noether(L, t, x, v, %);
\end{verbatim}
\begin{equation*}
0 = const
\end{equation*}
\normalsize
\end{explTEMAp}

\begin{explTEMAp}[Oscilador Harmónico com Amortecimento]
\label{ex:OscHarAmort}
Consideremos um oscilador harmónico com força de restituição $-kx$,
submerso num líquido de tal modo que o movimento da massa $m$
é amortecido por uma força proporcional à sua velocidade.
Recorrendo à segunda lei de Newton obtém-se, como equação de movimento,
a equação diferencial de Euler-Lagrange associada ao seguinte Lagrangeano
\cite[pp. 432--434]{LOG}:
\small
\begin{verbatim}
> L:=1/2 * (m*v^2-k*x^2)*exp((a/m)*t);
\end{verbatim}
\begin{equation*}
L := \frac{1}{2} \left( m{v}^{2}-k{x}^{2} \right) {e^{{\frac {at}{m}}}}
\end{equation*}
\normalsize
Para determinar um primeiro integral da equação de Euler-Lagrange,
encontramos os geradores $T$ e $X$ sob os quais a funcional integral
$J[x(\cdot)] = \int L dt$ é invariante:
\small
\begin{verbatim}
> Simetria(L, t, x, v);
\end{verbatim}
\begin{equation*}
 \left\{ T \left( t,x \right) ={\it \_C1},X \left( t,x \right) =
 - \frac{xa{\it \_C1}}{2m} \right\}
\end{equation*}
\normalsize
\small
\begin{verbatim}
> S:= subs(_C1=1,%);
\end{verbatim}
\begin{equation*}
 S:=  \left\{ X \left( t,x \right)
 =-{\frac {xa}{2m}},T \left( t,x \right) =1 \right\}
\end{equation*}
\normalsize
Pelo Teorema de Noether (Corolário~\ref{cor:TeoNoether}) obtemos o primeiro integral
$\left(L - \dot{x} \frac{\partial L}{\partial \dot{x}}\right)T$
$+ \frac{\partial L}{\partial \dot{x}} X$:
\small
\begin{verbatim}
> simplify(Noether(L, t, x, v, S));
\end{verbatim}
\begin{equation*}
-\frac{1}{2}\,{e^{{\frac {at}{m}}}} \left( x \left( t \right) a{\frac {d}{dt}}
x \left( t \right) +m \left( {\frac {d}{dt}}x \left( t \right)
 \right) ^{2}+k \left( x \left( t \right)  \right) ^{2} \right) ={\it const}
\end{equation*}
\normalsize
\end{explTEMAp}


\newsec{Procedimentos Maple}
\label{sec:ProcMaple}
\setcounter{obsTEMAp}{0}

Os procedimentos \emph{Simetria}, \emph{Noether} e \emph{EulerLagrange},
descritos e ilustrados nas secções anteriores, são agora definidos
usando o sistema de computação algébrica \texttt{Maple} (versão~9).

\begin{description}
\item[Simetria] determina as simetrias de um Lagrangeano
de várias variáveis dependentes e com derivadas de ordem superior,
de acordo com a secção~\ref{sec:SV}.
\item Devolve:
\begin{itemize}
\item[-] conjunto/vector de geradores infinitesimais
das transformações simétricas.
\end{itemize}
\item Forma de invocação:
\begin{itemize}
\item[-] Simetria(L, t, x, x1, x2, ..., xm)
\end{itemize}
\item Parâmetros:
\begin{itemize}
\item[L -] expressão do Lagrangeano;
\item[t -] nome da variável independente;
\item[x -] nome, lista de nomes ou vector de nomes das variáveis dependentes;
\item[xi -] nome, lista de nomes ou vector de nomes das derivadas de ordem i
das variáveis dependentes;
\end{itemize}
\end{description}
\small
\begin{verbatim}
Simetria:=proc(L::algebraic,t::name,x0::{name,list(name),'Vector[column]'
                       (name)},x1::{name,list(name),'Vector[column]'(name)})
  local n,m,xx,P,EqD,SysEqD,Sol,xi,Tdt,soma,V,r,i,j,k;
  if nargs<4 then print(`Nº de args insuficiente.`); return;
  elif not type([args[3..-1]],{'list'(name),'listlist'(name),
                                            'list'('Vector[column]'(name))})
    then print(`Erro na lista das var. depend. ou suas derivadas.`); return;
  end if;
  unassign('T'); unassign('X');
  xx:=convert(x0,'list')[]; n:=nops([xx]); m:=nargs-3;
  xi:=[seq(Vector(convert(args[i],'list')),i=3..m+3)];
  Tdt:=diff(T(t,xx),t)+Vector[row]([seq(diff(T(t,xx),i),i=xx)]).xi[2];
  if n>1 then P:=[Vector([seq(X[i](t,xx),i=1..n)])];
  else P:=[Vector([X(t,xx)])]; end if;
  for i from 1 to m do
    V:=Vector(n):
    for k from 0 to i-1 do
      V:=V+Matrix(n,(r,j)->diff(P[i][r], xi[k+1][j])).xi[k+2];
    end do:
    P:=[P[],map(diff,P[i],t)+V-xi[i+1]*Tdt]:
  end do:
  soma:=0:
  for i from 0 to m do
    soma:=soma+Vector[row]([seq(diff(L,r),r=convert(xi[i+1],list))]).P[i+1]:
  end do:
  EqD:=diff(L,t)*T(t,xx)+soma+L*Tdt;
  EqD:=collect(EqD,[seq(convert(xi[i+1],'list')[],i=1..m)],distributed);
  SysEqD:={coeffs(EqD,[seq(convert(xi[i+1],'list')[],i=1..m)])};
  Sol:=pdsolve(SysEqD,{T(t,xx)} union convert(P[1],'set'));
  if type(x0,'Vector') then
    return (subs(Sol,T(t,xx)),Vector(subs(Sol,P[1])));
  else return Sol; end if;
end proc:
\end{verbatim}
\normalsize

\begin{description}
\item[Noether] dados os geradores infinitesimais,
determina a Lei de Conservação de um Lagrangeano de várias
variáveis dependentes e com derivadas de ordem superior,
de acordo com a secção~\ref{sec:LC}.
\item Devolve:
\begin{itemize}
\item[-] lei de conservação.
\end{itemize}
\item Formas de invocação:
\begin{itemize}
\item[-] Noether(L, t, x, x1, x2, ..., xm, S)
\item[-] Noether(L, t, xs, x1s, x2s, ..., xms, T, X)
\end{itemize}
\item Parâmetros:
\begin{itemize}
\item[L -] expressão do Lagrangeano;
\item[t -] nome da variável independente;
\item[x -] nome ou lista de nomes das variáveis dependentes;
\item[xi -] nome ou lista de nomes das derivadas de ordem i
das variáveis dependentes;
\item[S -] conjunto de geradores infinitesimais das simetrias
(output do procedimento \textsl{Simetria});
\item[xs -] vector de nomes das variáveis dependentes;
\item[xis -] vector de nomes das derivadas de ordem i
das variáveis dependentes;
\item[T -] gerador da transformação para a variável independente (t);
\item[X -] vector com os geradores das transformações para as variáveis
dependentes (xs).
\end{itemize}
\end{description}
\small
\begin{verbatim}
Noether:=proc(L::algebraic,t::name,x0::{name,list(name),'Vector[column]'
                       (name)},x1::{name,list(name),'Vector[column]'(name)})
  local xx,n,m,P,psi,LC,xi,Tdt,V,r,i,j,k;
  if type(x0,'Vector') then m:=nargs-5; else m:=nargs-4; end if;
  if m<1 then print(`Nº de args insuficiente.`); return;
  elif not type([args[3..3+m]],{'list'(name),'listlist'(name),
                                       'list'('Vector[column]'(name))}) then
    print(`Erro na lista das var. depend. ou suas derivadas.`); return;
  elif (type(x0,'Vector') and not(type(args[-1],'Vector[column]') and
                    type(args[-2],algebraic))) or (not type(x0,'Vector') and
                                              not type(args[-1],'set')) then
    print(`Conj. de gerad. inválido.`); return;
  end if;
  xx:=convert(x0,'list')[];  n:=nops([xx]); unassign('T'); unassign('X');
  xi:=[seq(Vector(convert(args[i],'list')),i=3..m+3)];
  xi:=[xi[],seq(Vector([seq(x||i[k],k=1..n)]),i=m+1..2*m-1)];
  Tdt:=diff(T(t,xx),t)+Vector[row]([seq(diff(T(t,xx),i),i=xx)]).xi[2];
  if n>1 then P:=[Vector([seq(X[i](t,xx),i=1..n)])]:
  else P:=[Vector([X(t,xx)])]: end if:
  for i from 1 to (m-1) do
    V:=Vector(n):
    for k from 0 to i-1 do
      V:=V+Matrix(n,(r,j)->diff(P[i][r],xi[k+1][j])).xi[k+2] end do:
    P:=[P[],map(diff,P[i],t)+V-xi[i+1]*Tdt]:
  end do:
  psi:=[Vector[row]([seq(diff(L,i),i=convert(xi[m+1],'list'))])]:
  for i from m by -1 to 2 do
    V:=Vector[row](n):
    for k from 0 to 2*m-i do
      V:=V+LinearAlgebra[Transpose](xi[k+2]).Matrix(n,(j,r)->
                                        diff(psi[1][r],xi[k+1][j])); end do:
    psi:=[Vector[row]([seq(diff(L,i),i=convert(xi[i],'list'))])
                                               -map(diff,psi[1],t)-V,psi[]]:
  end do:
  LC:=sum('psi[i].P[i]','i'=1..m)+(L-sum('psi[i].xi[i+1]','i'=1..m))*T(t,xx)
                                                                     =const:
  if type(x0,'Vector') then
    LC:=eval(LC,[T(t,xx)=args[-2],seq(P[1][i]=args[-1][i],i=1..n)]);
  else LC:=eval(LC,args[-1]); end if;
  LC:=subs({map(i->i=i(t),[xx])[]},LC);
  LC:=subs({seq(seq(xi[k+1][i]=diff(xi[1][i](t),t$k),i=1..n),k=1..2*m-1)},
                                                                        LC);
  return LC;
end proc:

\end{verbatim}
\normalsize

\begin{description}
\item[EulerLagrange] constrói o sistema de equações
de Euler-Lagrange \eqref{eq:EL}, dado um Lagrangeano de várias variáveis
dependentes e com derivadas de ordem superior.
\item Devolve:
\begin{itemize}
\item[-] conjunto/vector de equações de Euler-Lagrange.
\end{itemize}
\item Forma de invocação:
\begin{itemize}
\item[-] EulerLagrange(L, t, x, x1, x2, ..., xm)
\end{itemize}
\item Parâmetros:
\begin{itemize}
\item[L -] expressão do Lagrangeano;
\item[t -] nome da variável independente;
\item[x -] nome, lista de nomes ou vector de nomes das variáveis dependentes;
\item[xi -] nome, lista de nomes ou vector de nomes das derivadas de ordem i
das variáveis dependentes;
\end{itemize}
\end{description}
\small
\begin{verbatim}
EulerLagrange:=proc(L::algebraic,t::name,x0::{name,list(name),
       'Vector[column]'(name)},x1::{name,list(name),'Vector[column]'(name)})
  local xx,n,m,Lxi,xi,V,EL,i,j,k;
  if nargs<4 then print(`Nº de args insuficiente.`); return;
  elif not type([args[3..-1]],{'list'(name),'listlist'(name),
                                            'list'('Vector[column]'(name))})
    then print(`Erro na lista das var. depend. ou suas derivadas.`); return;
  end if;
  xx:=convert(x0,'list')[]; n:=nops([xx]); m:=nargs-3;
  xi:=[seq(Vector(convert(args[i],'list')),i=3..m+3)];
  V:=[0$n];
  for i from 1 to m do
    Lxi:=[seq(diff(L,k),k=convert(xi[i+1],'list'))]:
    Lxi:=subs({map(k->k=k(t),[xx])[]},Lxi);
    Lxi:=subs({seq(seq(xi[k+1][j]=diff(xi[1][j](t),t$k),j=1..n),k=1..m)},
                                                                       Lxi);
    V:=V+(-1)^i*map(diff,Lxi,t$i);
  end do:
  EL:=[seq(diff(L,k),k=convert(xi[1],'list'))];
  EL:=subs({map(k->k=k(t),[xx])[]},EL);
  EL:=subs({seq(seq(xi[k+1][j]=diff(xi[1][j](t),t$k),j=1..n),k=1..m)},EL);
  EL:=EL+V;
  if type(x0,'Vector') then return convert(map(i->i=0,EL),'Vector[column]');
  elif type(x0,'list') then return convert(map(i->i=0,EL),'set');
  else return op(EL)=0; end if;
end proc:
\end{verbatim}
\normalsize


\newsec{Conclusão e Trabalho Futuro}
\label{sec:conclusao}
\setcounter{obsTEMAp}{0}

Os sistemas actuais de computação algébrica
colocam à nossa disposição ambientes de computação científica
extremamente sofisticados e poderosos. Disponibilizam
já muito conhecimento matemático e permitem estender
esse conhecimento por intermédio de linguagens de programação de
muito alto nível, expressivas e intuitivas, próximas da linguagem matemática.
Tais ambientes permitem a realização de uma miríade de cálculos matemáticos
simbólicos, com extrema eficiência, e a definição célere de
novas funcionalidades. Neste trabalho usámos o sistema de computação matemática
\texttt{Maple~9} para definir novas funções que permitem
a determinação automática de simetrias e leis de conservação
no cálculo das variações. Mostrámos depois, por meio de exemplos concretos,
como as novas funcionalidades são de grande utilidade prática.
Como trabalho futuro pretendemos estender
o nosso ``package Maple'' ao caso discreto e aos problemas mais
genéricos tratados pelo controlo óptimo.


\bigskip

Uma abordagem análoga à que fizemos pode ser efectuada
para problemas do cálculo das variações discretos no tempo.
No caso discreto, o problema fundamental do
cálculo das variações consiste em determinar uma sequência finita
$\mathbf{x}(k)\in \mathbb{R}^n$, $k=M,\ldots,M+N$, de  modo
a que a função de custo discreta
\begin{equation*}
J[\mathbf{x}(\cdot)] = \sum_{k=M}^{M+N-1}
L(k,\mathbf{x}(k),\mathbf{x}(k+1))
\end{equation*}
seja minimizada (ou maximizada). Neste tipo de problemas consideramos
um intervalo de tempo de $N$ períodos, com início num período fixo $M$,
em que $k\in\mathbb{Z}$, $k=M,\ldots,M+N-1$, é a variável discreta
tempo, e assumimos que o Lagrangeano $L$ é continuamente diferenciável
em $\{M,\ldots,M+N-1\}\times \mathbb{R}^n \times\mathbb{R}^n$.
À semelhança do que fizemos para o caso contínuo, é natural
considerarem-se problemas de optimização em que o Lagrangeano
envolve diferenças finitas de ordem superior,
\begin{equation}
\label{eq:FCDiscOrdSup}
J[\mathbf{x}(\cdot)] = \sum_{k=M}^{M+N-1}
L(k,\mathbf{x}(k),\mathbf{x}(k+1),\ldots,\mathbf{x}(k+m)),
\end{equation}
com $m\ge 1$, e considerando ainda o Lagrangeno continuamente diferenciável
relativamente a todos os seus argumentos.
Uma condição necessária para que $\mathbf{x}(k)$ seja uma extremal de
(\ref{eq:FCDiscOrdSup}) é dada pela equação de Euler-Lagrange discreta
\begin{equation}
\label{eq:EPDisc}
\sum_{j=0}^m\frac{\partial L}{\partial \mathbf{x}^j}
(k+m-j,\mathbf{x}(k+m-j),
\ldots,\mathbf{x}(k+2m-j))=\mathbf{0} \, .
\end{equation}
No caso $m=1$, a equação (\ref{eq:EPDisc}) reduz-se
à equação discreta de Euler-Lagrange
\begin{equation*}
\frac{\partial L}{\partial \mathbf{x}}(k+1,\mathbf{x}(k+1),
\mathbf{x}(k+2))
+\frac{\partial L}{\partial \mathbf{x}^1}
(k,\mathbf{x}(k),\mathbf{x}(k+1))=\mathbf{0}.
\end{equation*}
\begin{defTEMAp}
\label{defin:invarDisc}
Um problema discreto expresso pela função de custo (\ref{eq:FCDiscOrdSup})
diz-se invariante sob as transformações uni-paramétricas
$\mathbf{h}^s(k,\mathbf{x})$, com $\mathbf{h}^0(k,\mathbf{x})=
\mathbf{x}$,
se, para todo o $s$ suficientemente pequeno e qualquer que seja o $k$,
\begin{equation}
\label{eq:invarDisc}
L(k,\mathbf{x}(k),\mathbf{x}(k+1),\ldots,\mathbf{x}(k+m))
= L(k,\mathbf{h}^s(k,\mathbf{x}(k)),\ldots,\mathbf{h}^s(k+m,\mathbf{x}(k+m))).
\end{equation}
\end{defTEMAp}
\noindent A partir desta definição conseguimos estabelecer uma condição necessária e
suficiente de invariância.
\begin{teoTEMA}
\label{teor:cnsDisc}
A função de custo (\ref{eq:FCDiscOrdSup})
é invariante sob as transformações uni-paramétricas
$\mathbf{h}^s(k,\mathbf{x})$, com geradores infinitesimais
$\mathbf{X}(k,\mathbf{x})=\left.\frac{\partial}{\partial{s}}
\mathbf{h}^s(k,\mathbf{x})\right|_{s=0}$,
se, e apenas se,
\begin{eqnarray}
\label{eq:cnsDisc}
\sum_{j=0}^{m}\frac{\partial L}{\partial \mathbf{x}^{j}}\cdot
\mathbf{X}(k+j,\mathbf{x}(k+j))=0,
\end{eqnarray}
onde o Lagrangeano discreto $L$ é avaliado em
$(k,\mathbf{x}(k),\mathbf{x}(k+1),\ldots,\mathbf{x}(k+m))$.
\end{teoTEMA}
\begin{obsTEMAp}
Quando $m=1$, a equação (\ref{eq:cnsDisc}) reduz-se a
\begin{equation*}
\frac{\partial L}{\partial \mathbf{x}}\cdot
\mathbf{X}(k,\mathbf{x}(k))
+\frac{\partial L}{\partial \mathbf{x}^1}
\cdot\mathbf{X}(k+1,\mathbf{x}(k+1))=0.
\end{equation*}
\end{obsTEMAp}
\noindent \emph{Demonstração (Teorema~\ref{teor:cnsDisc})}.
Da condição de invariância (\ref{eq:invarDisc}) podemos escrever:
\[
\frac{\textrm{d}}{\textrm{d}s}L(k,\mathbf{h}^s(k,\mathbf{x}(k)),
\ldots,\mathbf{h}^s(k+m,\mathbf{x}(k+m)))=0 \, .
\]
Derivando, obtemos, para $s = 0$, a igualdade \eqref{eq:cnsDisc}:
\begin{displaymath}
\sum_{j=0}^{m}\frac{\partial L}{\partial \mathbf{x}^{j}}\cdot
\left.\frac{\partial}{\partial{s}}\mathbf{h}^s(k+j,\mathbf{x}(k+j))
\right|_{s=0}=0
\end{displaymath}
com $L$ avaliado em
$(k,\mathbf{x}(k),\mathbf{x}(k+1),\ldots,\mathbf{x}(k+m))$.
\begin{flushright}
$\square$
\end{flushright}
É possível deduzir um teorema análogo ao
teorema de Noether, para o caso discreto, que permite estabelecer
leis de conservação discretas (cf. \cite{delfimCAO}).
\begin{teoTEMA}[Teorema de Noether Discreto]
\label{teor:lcDisc}
Se a função de custo (\ref{eq:FCDiscOrdSup}) é invariante,
no sentido da definição \ref{defin:invarDisc},
sob as transformações uni-paramétricas
$\mathbf{h}^s(k,\mathbf{x})$, com geradores infinitesimais
$\mathbf{X}(k,\mathbf{x})=\left.\frac{\partial}{\partial{s}}
\mathbf{h}^s(k,\mathbf{x})\right|_{s=0}$,
então todas as soluções $\mathbf{x}(k)$
da equação de Euler-Lagrange (\ref{eq:EPDisc}) satisfazem
\begin{eqnarray}
\label{eq:LCDiscOrdSup}
\sum_{j=0}^{m-1}\mathbf{\Psi}^j(k)\cdot
\mathbf{X}(k+j,\mathbf{x}(k+j))=const\quad
\end{eqnarray}
onde
\begin{eqnarray}
\mathbf{\Psi}^0(k)&=&\frac{\partial L}{\partial \mathbf{x}}
(k,\mathbf{x}(k),\ldots,\mathbf{x}(k+m))\nonumber\\
\mathbf{\Psi}^{j}(k)&=&\mathbf{\Psi}^{j-1}(k+1)
+\frac{\partial L}{\partial \mathbf{x}^{j}}
(k,\mathbf{x}(k),\ldots,\mathbf{x}(k+m))\textrm{,}\quad
para \; j=1,2,\ldots,m-1\,.\nonumber
\end{eqnarray}
\end{teoTEMA}
\begin{obsTEMAp}
Caso $m=1$, a lei de conservação (\ref{eq:LCDiscOrdSup}) reduz-se a
\begin{equation*}
\frac{\partial L}{\partial \mathbf{x}}(k,\mathbf{x}(k),
\mathbf{x}(k+1))\cdot
\mathbf{X}(k,\mathbf{x}(k))=const.
\end{equation*}
\end{obsTEMAp}
\noindent Para problemas discretos no tempo, facilmente se definem
em Maple procedimentos que implementem as condições
(\ref{eq:EPDisc}), \eqref{eq:cnsDisc} e (\ref{eq:LCDiscOrdSup}).


\bigskip

O controlo óptimo pode ser encarado como uma extensão natural
do cálculo das variações \cite{CVFAwOC,Leitao}.
A resolução de problemas do controlo óptimo
passa normalmente pela aplicação do Princípio do Máximo
de Pontryagin, que constitui uma generalização
das condições clássicas de Euler-Lagrange
e de Weierstrass do cálculo das variações.
Em termos práticos, algorítmicos, temos
de resolver um sistema de equações diferenciais, chamado de sistema
Hamiltoniano, depois de eliminados os parâmetros
de controlo por intermédio da \emph{condição de máximo}.
Tal como no cálculo das variações, também no controlo
óptimo as equações diferenciais ordinárias obtidas são em geral
não-lineares e de difícil resolução (podem,
inclusive, não ser integráveis). As leis
de conservação são usadas para simplificar essas equações e,
mais uma vez, a questão que se coloca é saber
como determinar tais quantidades preservadas.
Resulta que os resultados clássicos de Emmy Noether
podem ser generalizados ao contexto do controlo óptimo,
reduzindo o problema ao da descoberta de grupos de transformações
uni-paramétricas que deixem o problema de controlo óptimo
invariante \cite{delfimEJC}. Seria de grande utilidade
prática dispor de facilidades computacionais
que permitissem identificar as simetrias dos problemas de controlo óptimo
\cite{GrizzleMarcus,alik}.

\bigskip

Serão estas algumas das direcções do nosso trabalho futuro.


\end{document}